\documentclass[preprint,12pt]{elsarticle}

\usepackage{amssymb}
\usepackage{color, amsmath,amssymb, amsfonts, amstext,amsthm, latexsym}

\usepackage{epsfig, graphicx, graphics}
\usepackage{longtable}

\newcommand{\p}{\partial}

\renewcommand{\phi}{\varphi}

\newcommand{\R}{{\mathbb R}}

\usepackage{amssymb}

\journal{\emph{Acta Mathematica Scientia }(Revised version)}

\begin{document}

\begin{frontmatter}



\title{Elementary bifurcations for a simple dynamical system under non-Gaussian L\'evy noises \footnote{Corresponding author: Huiqin Chen,E-mail:chenhuiqin111@yahoo.com.cn.
This work was partly supported
by  the NSFC grants 10971225, 11171125, 91130003 and 11028102,  the NSFH  2011CDB289, HPDEP  20114503 and 2011B400,
the Cheung Kong Scholars Program and    the
Fundamental Research Funds for the Central Universities, HUST
2010ZD037.}}


\author{ Huiqin Chen$^1$$^,$ $^2$,  Jinqiao Duan $^{3}$ and  Chengjian Zhang$^1$}

\address{\emph{1. School of Mathematics and Statistics,  Huazhong University of Science and Technology,
  Wuhan 430074, China }\emph{2. School of Mathematics and Computer Science,  Jianghan University ,
  Wuhan 430056, China } E-mail: chenhuiqin111@yahoo.com.cn    \, cjzhang@mail.hust.edu.cn\\
  \emph{3. Department of Applied Mathematics, Illinois Institute of Technology,
  Chicago, IL 60616, USA} E-mail: duan@iit.edu  }

\begin{abstract}
Nonlinear dynamical systems are sometimes under the influence of random fluctuations. It is desirable to examine possible bifurcations for stochastic dynamical systems when a parameter varies.

A computational analysis is   conducted to investigate
 bifurcations of a simple dynamical system under non-Gaussian $\alpha-$stable L\'evy motions, by examining the  changes in   stationary probability density
 functions for the solution orbits of this stochastic system.  The stationary probability density
 functions are obtained by numerically solving a nonlocal Fokker-Planck equation. This allows numerically investigating phenomenological bifurcation,
 or P-bifurcation, for stochastic differential equations with non-Gaussian L\'evy noises.
\end{abstract}


\end{frontmatter}


\section{Motivation}  \label{intro}
The dynamical behaviors for a dynamical system  depending on a parameter may change when this parameter   varies. This so called bifurcation phenomenon has been observed in many deterministic systems \cite{GH}. It also occurs in stochastic systems with Gaussian noises \cite{Arnold}.
It is desirable to examine possible bifurcation phenomena for stochastic systems with non-Gaussian noises.

L\'evy motions $L_t$ are a class of stochastic processes that have independent and stationary increments. They are usually non-Gaussian processes.    The well-known Brownian motion $B_t$ is a special case  which has \emph{additional} properties:  (i) Almost every
sample path of the Brownian motion     is  continuous in time in the
usual sense and (ii) Brownian motion's increments have Gaussian
distribution. Random fluctuations in  complex systems in engineering and science are often  non-Gaussian.
For instance, it has been argued that diffusion by
geophysical turbulence \cite{Shlesinger} corresponds, loosely
speaking, to a series of  ``pauses", when the particle is trapped by
a coherent structure, and ``flights" or ``jumps" or other extreme
events, when the particle moves in the jet flow. Paleoclimatic data
\cite{Dit} also indicates such irregular processes.

  SDEs  perturbed by   non-Gaussian L\'evy noises
have attracted much attention  recently \cite{Apple,  Schertzer}.
    SDEs perturbed by L\'evy motion   generate stochastic
flows \cite{Kunita2004, Apple}, or   random dynamical
systems  (cocycles), under certain conditions.

\medskip

Let us   consider a deterministic   differential equation
perturbed   by a non-Gaussian L\'evy motion, i.e., consider a stochastic differential equation (SDE)
\begin{equation} \label{sde}
dX_{t}  =  f (b, X_{t})  dt + \epsilon d L_t^{\alpha},
\end{equation}
 where $b \in \R$,   $ \epsilon>0$  and $\alpha \in (0, 2)$ are real parameters;  and $L_t^{\alpha}$ is a $\alpha-$stable symmetric  L\'evy motion    defined in a probability space $(\Omega, \mathcal{F}, \mathbb{P})$.
In this paper, we consider a numerical approach for understanding how the dynamic behaviors change when parameters vary, for a special case $f=b X_t-X_t^3$. Note that  $\dot{x} = bx-x^3$ is a primary dynamical model exhibiting  the deterministic  pitchfork bifurcation \cite{GH}.

\vskip 12pt

In section 2, we   briefly review some basic concepts
for L\'evy motions. In section 3, we first discuss the nonlocal Fokker-Planck equations for SDEs with L\'evy motions, then    present a numerical approach in computing stationary probability densities for the solution processes for \eqref{sde} above,
 and further discuss how the  stationary probability densities change when the parameters vary (i.e., phenomenological bifurcation or P-bifurcation).

\section{L\'evy motions and generators}    \label{motion}

Let us briefly review basic facts about L\'evy motions.
A scalar L\'evy motion is characterized by a drift parameter
$\theta$, a variance (or diffusion) parameter $d \geq 0$ and   a non-negative Borel
measure $\nu$, defined on $(\R, \mathcal{B}(\R))$ and concentrated
on $\R \setminus\{0\}$, which satisfies
\begin{equation} \label{levycondition}
  \int_{\R \setminus\{0\} } (y^2 \wedge 1) \; \nu(dy) < \infty,
\end{equation}
or equivalently
\begin{equation}
  \int_{\R \setminus\{0\} } \frac{y^2}{1+y^2}\; \nu(dy) < \infty.
\end{equation}
This measure $\nu$ is the so called    L\'evy jump measure of
L\'evy motion $L(t)$. We also call $(\theta, d, \nu)$ the
\emph{generating triplet}.

Let $L_t$ be a L\'evy process with the generating triplet
$(\theta, d, \nu)$. It is known that a scalar L\'evy
motion is completely determined by the L\'evy-Khintchine formula
(See \cite{Apple, Sato-99, PZ}). This says that for any
one-dimensional L\'evy process $L_{t}$, there exists a $\theta \in
\R$, $d>0$ and a measure $\nu$ such that its characteristic function is
\begin{equation}
Ee^{i \lambda L_{t}}=\exp \{i \theta\lambda t - d t
\frac{\lambda^{2}}{2} +t \int_{\R \setminus \{0\} } (e^{ i \lambda
y}-1 -i \lambda y I_{\{|y| <1 \}} ) \nu(dy)\},
\end{equation}
where $I_S$ is the indicator function of the set $S$, defined as follows:
$$
I_S(y) =
    \begin{cases}
        1,   &\text{if $y \in S $;}\\
        0,    &\text{if $y \notin S$.}
    \end{cases}
$$

The generator $\tilde{A}$ of the process $L_{t}$ is the same as the
infinitesimal generator since L\'evy process has independent and
stationary increments. Hence  $\tilde{A}$ is defined as $\tilde{A}  \phi = \lim_{t
\downarrow 0} \frac{P_{t} \phi -\phi}{t}$ where $P_{t} \phi(x)=
E_{x} \phi(L_{t})$ and $\phi$ is any function belonging to the
domain of the operator $\tilde{A}$. Recall the generator $\tilde{A}$ for  $L_{t}$ is
(See \cite{Apple, PZ})
\begin{equation} \label{A}
\tilde{A}  \phi(x) =  \theta \phi'(x) + \frac{d}{2} \phi''(x) + \int_{\R
\setminus\{0\}} [\phi(x+ y)-\phi(x) -  I_{\{|y|<1\}} \; y \phi'(x)
] \; \nu(dy).
\end{equation}

 \medskip

  In this paper, we consider a special L\'evy process, the symmetric $\alpha-$stable L\'evy motion $L_t^{\alpha}$, with drift $\theta=0$, diffusion $d=0$ and the jump measure $\nu_{\alpha}(dy)= \frac{dy}{|y|^{1+\alpha}}$. The corresponding generator is
 $A_{\alpha}  \phi(x) =    \int_{\R
\setminus\{0\}} [\phi(x+ y)-\phi(x) -  I_{\{|y|<1\}} \; y \phi'(x)
]\nu_{\alpha}(dy) $.

In the  next section, we consider bifurcation of the    equation \eqref{sde} when $b$ and $\alpha$ vary, by numerically investigating stationary probability density function for the solution of \eqref{sde}. We take the drift coefficient $f(b, x) = bx-x^3$, corresponding to the well-known pitchfork bifurcation in the deterministic case (when $\epsilon=0$). Note that simulations for solution paths were conduced for SDEs with $\alpha-$stable L\'evy noises in \cite{JW}, while we examine bifurcation phenomena by computing stationary probability density functions for solutions in the present paper.


\section{Bifurcation under additive L\'evy noises }\label{additive}

We now consider possible bifurcations for the SDE
\begin{equation} \label{sde2}
dX_{t}  =  f (b, X_{t})  dt + \epsilon d L_t^{\alpha}.
\end{equation}

The generator $A$ for the solution process   $X_{t} $  in \eqref{sde2} is
\begin{eqnarray} \label{AA0}
A  \phi = f(b, x) \phi^{'}(x)
 +\epsilon \int_{\R \setminus\{0\}} [\phi(x+  y)-\phi(x)-  I_{\{|y|<1\}} \; y \phi'(x) ] \; \nu_{\alpha}(dy).
\end{eqnarray}
 The adjoint operator for $A$ in the Hibert space $L^2(\R)$, with the usual scalar product, is then \cite{wu}
 \begin{eqnarray} \label{A*}
A^*  \phi = - [f(b, x) \phi(x)]^{'}
 + \epsilon\int_{\R \setminus\{0\}} [\phi(x+  y)-\phi(x) - I_{\{|y|<1\}} \; y \phi'(x)] \; \nu_{\alpha}(dy).
\end{eqnarray}
Consequently, the Fokker-Planck equation for the probability density function $p(x, t)$ for the solution
 process   $X_{t} $  in \eqref{sde} is
 \begin{eqnarray} \label{FPE}
\p_t p &=& - \p_x [f(b, x) p(x, t)]   \nonumber \\
 &+& \epsilon\int_{\R \setminus\{0\}} [p(x+  y, t)-p(x, t) - I_{\{|y|<1\}}\; y \; \p_x p(x, t) ]] \; \nu_{\alpha}(dy).
\end{eqnarray}

The   stationary solutions $p(x)$ of the above Fokker-Planck equation define  some invariant measures for the equation \eqref{sde}:
for real parameter $b$ and L\'evy parameter $\alpha \in (0, 2)$. Namely, the stationary probability density function $p(x)$ satisfies
\begin{eqnarray} \label{steady}
  - [f(b, x) p(x)]^{'}
  +  \epsilon\int_{\R \setminus\{0\}} [p(x+  y)-p(x) - I_{\{|y|<1\}}y p'(x)] \; \frac{dy}{|y|^{1+\alpha}}= 0, \\
   p(x) \geq 0, \;  \int_{\R} p(x) dx =1.
\end{eqnarray}

\subsection{Deterministic  pitchfork bifurcation}
   The differential equation $\dot x_{t}  = bx_t -x_t^3$ has a stable   fixed point  at $x=0$ for $ b <0$, and two additional stable   fixed points  at $ x=\pm\sqrt{b}$ for  $ b >0$ (See \cite{GH}). It undergoes a  so-called pitchfork bifurcation at $b=0$.

\subsection{Bifurcation under Brownian motion} \label{bifn-BM}

 We  first recall a  bifurcation under Brownian motion (See \cite{CF} or \cite[Page 475]{Arnold}), i.e., in the  case when the L\'evy motion in the above equation \eqref{sde2} is replaced by a Brownian motion $W_t$.  For $dX_{t}  = (bX_t -X_t^3) dt + \sigma dW_t,b\in R,\sigma \neq 0$, there exists only a unique stationary measure with density   $p_{b,\sigma}=N_{b,\sigma} \exp[\frac{1}{\sigma^2}(bx^2-\frac{x^4}{2})]$, where $N_{b,\sigma}$ is a normalization constant.
  For any give noise intensity $\sigma\neq0$,  the density is unimodal for $b \leq 0$, but    bimodal for $ b>0$
  (and the plateau for $p(x)$ occurs at $x_1=\sqrt{b}$ and $x_2=-\sqrt{b}$).  Hence the family ( ${p_{(b,\sigma)}}_{b\in R} $)
 undergoes a  bifurcation at $b=0$ for each $\sigma \neq 0$. This is a kind of phenomenological bifurcation or P-bifurcation \cite{Arnold},  in which a Brownian motion leads to a different bifurcation than its deterministic counterpart.

\subsection{Bifurcation under $\alpha-$stable L\'evy motion}\label{bifn-BMm}

 Now we consider the  bifurcation under $\alpha-$stable L\'evy motion. For
 $dX_{t}  = (bX_t -X_t^3) dt + \epsilon d L_t^{\alpha},b\in R,\alpha\in (0, 2)$,   a  stationary measure has density $p(x)=p_{b,\epsilon,\alpha}(x)$
    satisfying the following steady Fokker-Planck equation
\begin{eqnarray} \label{steady3}
 - [(bx-x^3)p(x)]^{'}
 +  \epsilon\int_{\R \setminus\{0\}} [p(x+  y)-p(x) -I_{\{|y|<1\}}y p'(x)] \; \frac{dy}{|y|^{1+\alpha}} = 0, \\
  p(x) \geq 0, \; \int_{\R} p(x) dx =1.
\end{eqnarray}

Unlike the Brownian case in     \S \ref{bifn-BM} above, we do not have the exact analytical solution         for the equation  \eqref{steady3}. In order to detect possible bifurcations, we instead  numerically simulate this  integro-differential equation on the interval
$(-l, l)$, with $l>0$ large enough and with the homogeneous Dirichlet boundary conditions.
 This  integro-differential equation contains   both an differential  part $ - [f(b, x) p(x)]^{'}$ and an integral part $\int_{\R \setminus\{0\}} [p(x+  y)-p(x) +I_{\{|y|<1\}}y p'(x)] \; \frac{dy}{|y|^{1+\alpha}} $.
We use a finite difference scheme on the differential part and the trapezoid rule in the integral  part \cite{Otto}; see also \cite{ChenLevy}.



We conduct  numerical simulations for various $b \in \mathbb{R}, \alpha \in (0, 2)$ and $\epsilon \in (0, 1)$.
Although $b$ and $\epsilon$ may be any real number in our numerical approach, here we
limit them $b$ to be in a bounded interval  in this paper.
All figures are in color in the online   version of this paper. Different colors are used to distinguish  cases with various parameter values.
In the following   only some  selected figures are shown to illustrate our results.


\subsubsection{Varying the   parameters $b$ and $\alpha$ } \label{ba}

\begin{figure}
 \begin{center}
\scalebox{0.9 }{\includegraphics{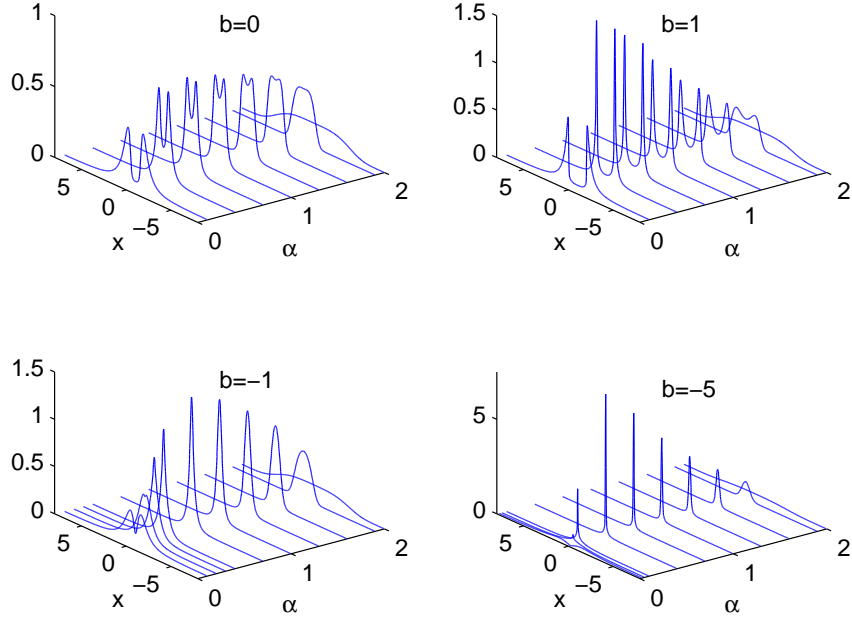}}
\caption{Stationary probability density   $p(x)$ for  $\epsilon=0.1$  }  \label{plot1}
 \end{center}
 \end{figure}

\begin{figure}
 \begin{center}
\scalebox{0.9}{\includegraphics{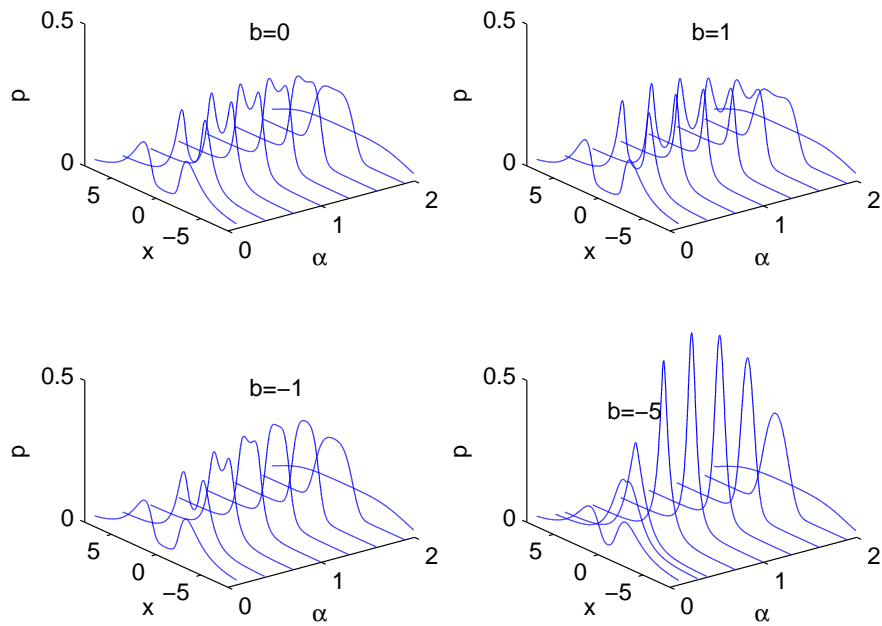}}
\caption{Stationary probability density   $p(x)$ for $\epsilon=0.9$  }  \label{plot2}
 \end{center}
 \end{figure}

Figures \ref{plot1} and \ref{plot2} show the stationary probability density function $ p(x)$    for $\epsilon=0.1$ and $\epsilon=0.9$, respectively.     Here we only show  several   cases for $ b=-5, -1, 0, 1$  and $\alpha=0.1, 0.4 , 0.7, 1.0, 1.3, 1.6, 1.9, 1.999$, as   examples.

The probability density function $p(x)$ evolves from bimodal  to unimodal, and  then further changes to the flatter kurtosis shape for every  fixed parameter $b$ in interval$ (-10, 0)$, as $\alpha$ value increases. The bifurcation occurs only for $b<0$, at some bifurcation value $\alpha*=\alpha(b, \epsilon)$, in our computational range. Note also that the  $p(x)$ is  bimodal for $ b\geq0$. This phenomenon is more evident when $\epsilon$ is larger in the range $(0, 1)$; see Figure 2.

For fixed $b$ value, the probability density function $p(x)$ evolves from  lower kurtosis  to the higher one, and then changes to lower again,   as $\alpha$ value increases.



\subsubsection{Varying the   parameter  $b$  } \label{b}

\begin{figure}
 \begin{center}
\scalebox{0.7}{\includegraphics{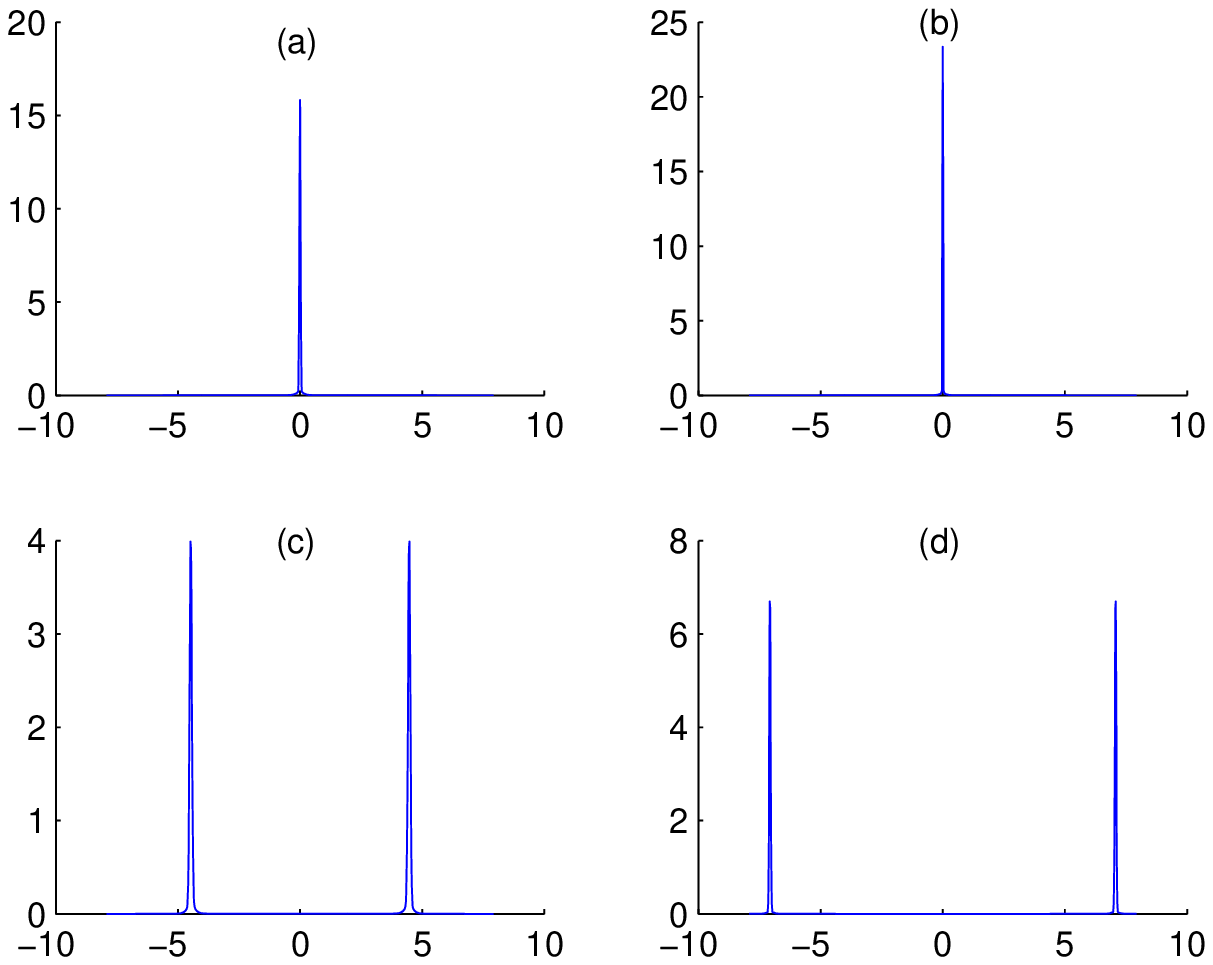}}
\caption{Stationary probability density function $p(x)$ for  $\epsilon=0.1$.  (a) Corresponding to $\alpha=1$ and $b=-30$; (b) Corresponding to $\alpha=1$ and $b=-50$; (c) Corresponding to $\alpha=1.5$ and $b=20$; (d) Corresponding to $\alpha=1.5$ and $b=50$}  \label{compare1}
 \end{center}
 \end{figure}

When $\alpha$ is approximately within  in the interval $(0.4, 1.6)$, the stationary density $p(x)$ becomes
very spiky, and this is more evident when the magnitude of $b$ is large; see Figure \ref{compare1}.

\subsubsection{Varying the   parameter  $\alpha$  } \label{b1}

\begin{figure}
 \begin{center}
\scalebox{0.7}{\includegraphics{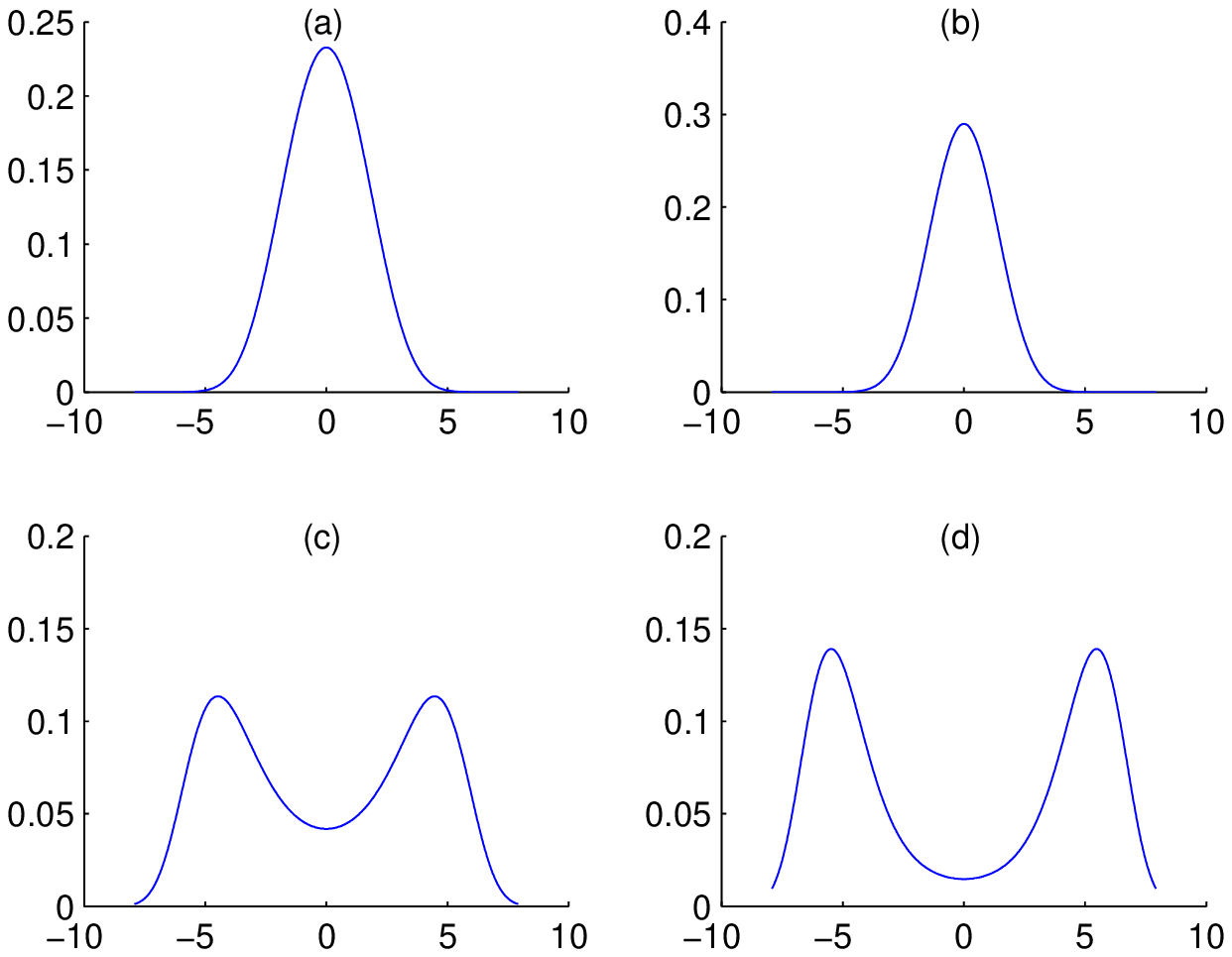}}
\caption{Stationary probability density function $p(x)$ for  $\epsilon=0.1$ and $\alpha=1.999$. (a) Corresponding to $b=-30$; (b) Corresponding to $b=-50$; (c) Corresponding to   $b=20$; (d) Corresponding to   $b=30$} \label{compare2}
 \end{center}
 \end{figure}

When $\alpha$ is approximately within  in the interval $(1.9, 2)$, the stationary density $p(x)$ becomes
very  flatter, and this is more evident when the value of $\alpha$ is close to 2; see Figure \ref{compare2}(subfigures (c) and (d) in both figures).

\subsubsection{Impact of noise intensity $\epsilon$   }\label{additive4}



\begin{figure}
 \begin{center}
\scalebox{0.7}{\includegraphics{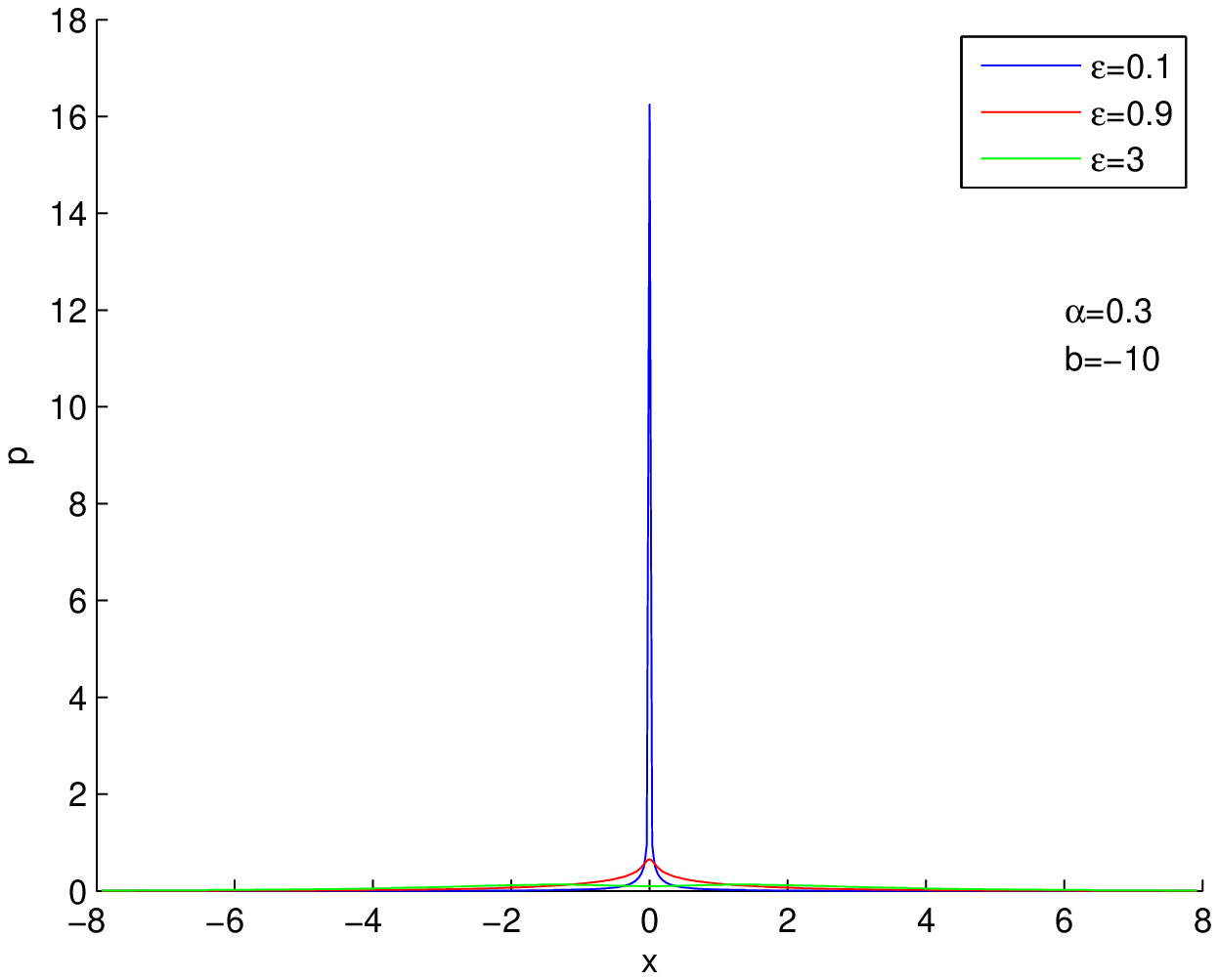}}
\caption{Stationary probability density   $p(x)$: $\alpha=0.3, \; b=-10$ for  $\epsilon=0.1$, $\epsilon=0.9$ and $\epsilon=3$  }  \label{good1}
 \end{center}
 \end{figure}

\begin{figure}
 \begin{center}
\scalebox{0.7}{\includegraphics{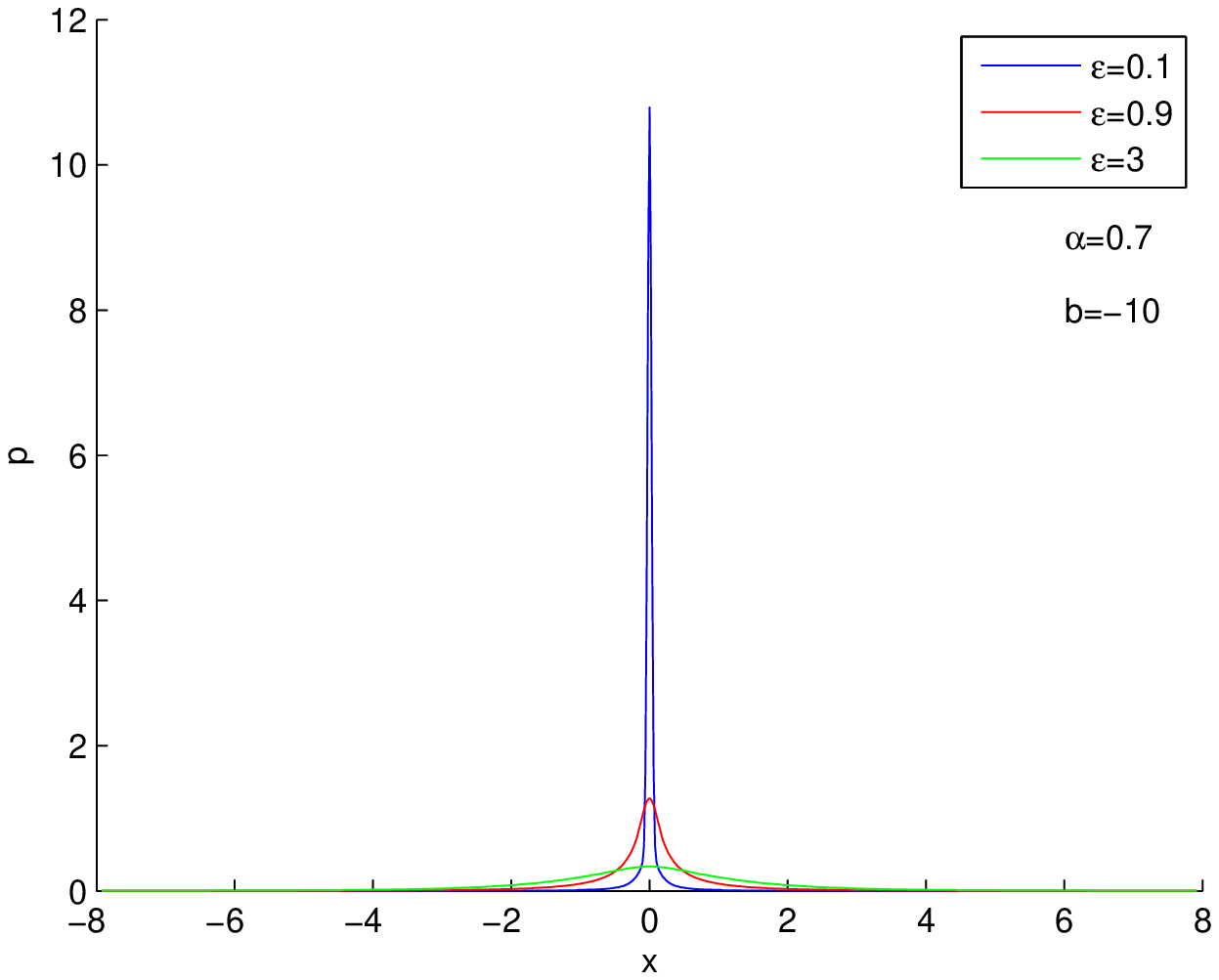}}
\caption{Stationary probability density   $p(x)$: $\alpha=0.7, \; b=-10$  for  $\epsilon=0.1$, $\epsilon=0.9$ and $\epsilon=3$  }  \label{good2}
 \end{center}
 \end{figure}

\begin{figure}
 \begin{center}
\scalebox{0.7}{\includegraphics{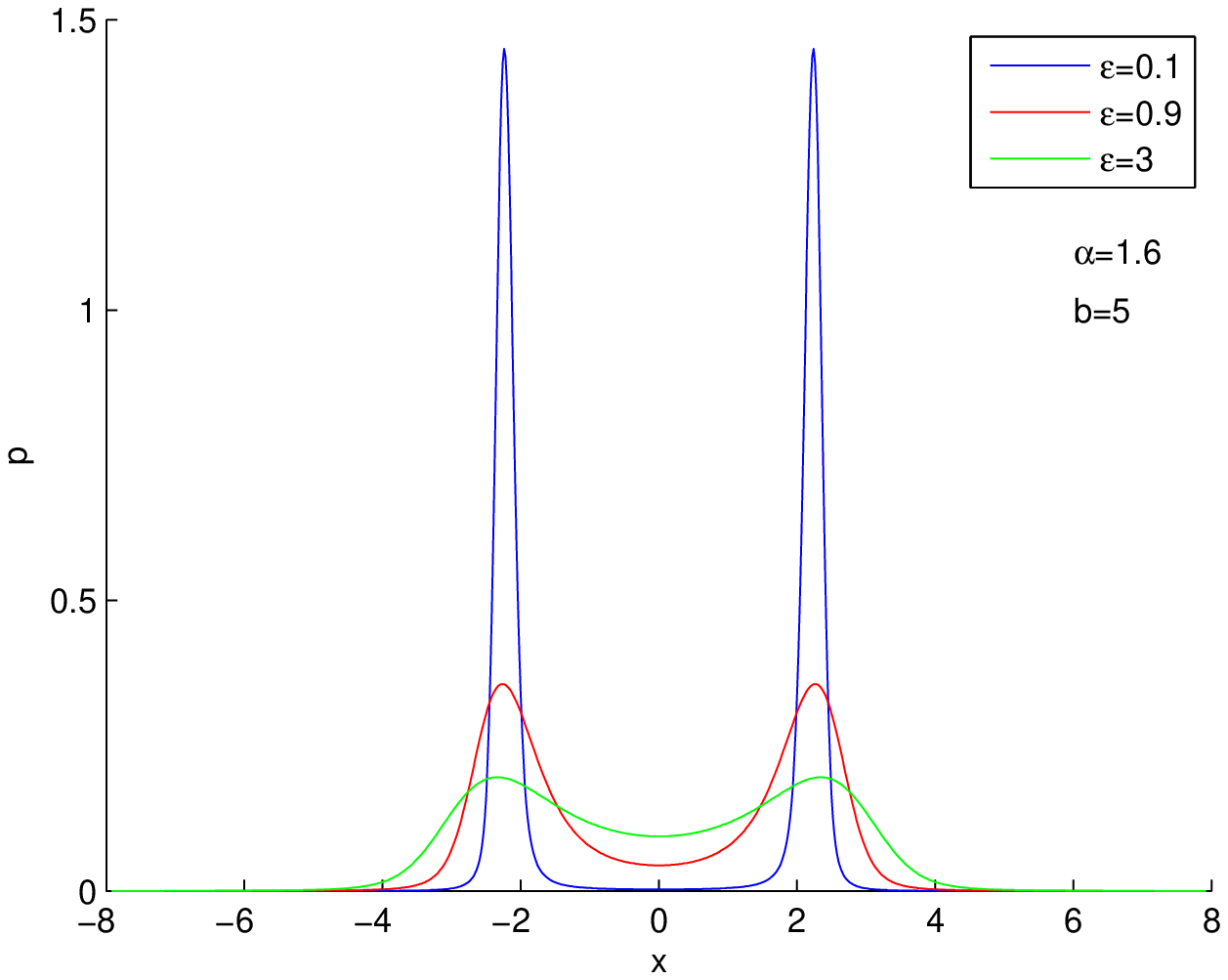}}
\caption{Stationary probability density   $p(x)$: $\alpha=1.6, \; b=5$  for  $\epsilon=0.1$, $\epsilon=0.9$ and $\epsilon=3$    }   \label{good3}
 \end{center}
 \end{figure}

\begin{figure}
 \begin{center}
\scalebox{0.7}{\includegraphics{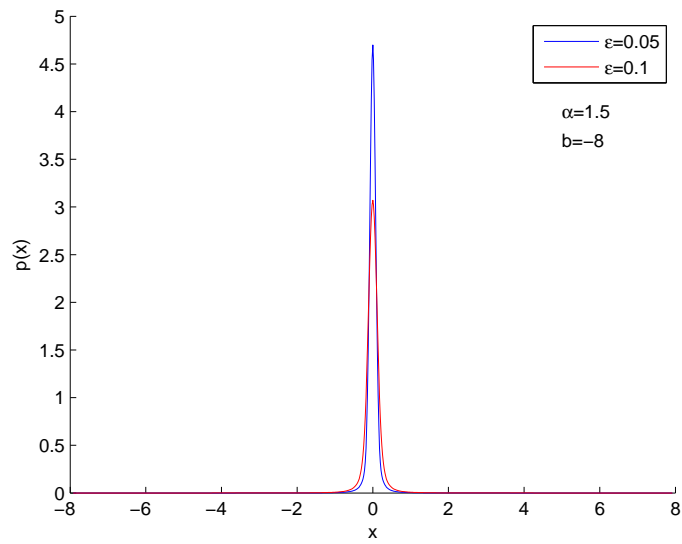}}
\caption{Stationary probability density   $p(x)$: $\alpha=1.5, \; b=-8$  for $\epsilon=0.05$ and $\epsilon=0.1$ } \label{good4}
 \end{center}
 \end{figure}

 When the positive noise intensity  $\epsilon$ is increased, we observe that the stationary density $p(x)$ becomes flatter (or less spiky) for fixed parameters $b$ and $\alpha$; see Figures \ref{good1}-\ref{good4}.
 Thus the nonlocal integral term, due to L\'evy jumps in the random forcing, in the steady Fokker-Planck equation \eqref{steady3} has a certain damping or diffusive effect.

\subsubsection{A remark:  Brownian motion vs. $\alpha-$stable L\'evy motion}


When $\alpha=2$, the corresponding L\'evy motion is the Brownian motion.
We take $\alpha$ close to $2$, the stationary density $p(x)$ are computed for various $b$ values; see Figures \ref{compare1} and \ref{compare2} (subfigures (c) and (d) in both figures), and also Figure \ref{good3}. Here we also observe that $p(x)$ is bimodal for $b>0$, similar to the Brownian motion case in \S
\ref{bifn-BM}.


\bigskip

{\bf Acknowledgements.} The authors would like to thank Xiaofan Li for helpful discussions about numerical schemes used in this research.






\bibliographystyle{elsarticle-num}


\end{document}